\documentclass[11pt]{amsart}
\usepackage{amsmath,amssymb,mathrsfs,pstricks}
\usepackage{tikz-cd} 
\usepackage[colorlinks=true,linkcolor=black,citecolor=black,urlcolor=black]{hyperref}

\usepackage{pxfonts}
\usepackage{setspace}
\usepackage[shortlabels]{enumitem}

\psset{unit=1pt}
\psset{arrowsize=4pt 1}
\psset{linewidth=.5pt}

\newcommand{\define}{\textit}

\newcommand{\excise}[1]{}

\renewcommand{\phi}{\varphi}

\newcommand{\PP}{\mathbb{P}}

\newcommand{\ZZ}{\mathbb{Z}}

\newcommand{\shfF}{\mathscr{F}}

\newcommand{\catA}{\mathcal{A}}

\newcommand{\catC}{\mathcal{C}}

\newcommand{\catE}{\mathcal{E}}
\newcommand{\catCoh}{\mathbf{Coh}}
\newcommand{\catM}{\mathcal{M}}

\newcommand{\OO}{\mathcal{O}}

\newcommand{\bn}{\mathbf{n}}
\newcommand{\bm}{\mathbf{m}}

\DeclareMathOperator{\im}{im}

\newtheorem{theorem}{Theorem}
\newtheorem{lemma}[theorem]{Lemma}

\theoremstyle{definition}

\begin{document}

\title{Gillet descent for connective K-theory}
\author{David Anderson}

\email{anderson.2804@math.osu.edu}

\date{July 21, 2020}
\thanks{DA partially supported by NSF CAREER DMS-1945212.}

\begin{abstract}
Using Gillet's technique of projective envelopes, we prove a homological descent theorem for the connective K-homology of schemes.
\end{abstract}

\maketitle

The aim of this note is to describe an exact descent sequence in connective K-theory.  This is a variation on a technique of Gillet, who used it to construct proper pushforwards for K-theory of general schemes, by bootstrapping from Quillen's construction for quasi-projective schemes.  Likewise, by using the present descent sequence, one can extend some results stated in \cite{cai,dl} for quasi-projective schemes to general schemes, at least for the part of  $CK_*(X)$ which specializes to $K_\circ(X)$ at $\beta=-1$.

The exposition will be terse: the arguments are all based on ones in the literature, and we refer especially to \cite{gillet} for the details.  See also \cite[Appendix~A]{ap} for a digest.

All schemes are separated and of finite type over some field.  Following Cai \cite{cai} (see also \cite[Appendix~A]{a}), let $\mathcal{M}_i(X) \subseteq \catCoh(X)$ be the full subcategory of sheaves whose support has dimension $\leq i$.  The \define{connective K-theory} groups of $X$ are defined as
\[
  CK_i(X) := \im\big( K_\circ(\catM_i(X)) \to K_\circ(\catM_{i+1}(X)) \big).
\]
More generally, Cai defines
\[
  CK_{i,q-i}(X) := \im\big((K_q( \catM_i(X)) \to K_q(\catM_{i+1}(X)) \big),
\]
although for $q>0$ his higher K-groups diverge from those of Dai and Levine \cite{dl}.  As mentioned above, here we focus only on the case $q=0$, where the notions from \cite{cai} and \cite{dl} agree.  

A proper morphism $f\colon Z \to X$ is an \define{envelope} if every subvariety of $X$ is the birational image of some subvariety of $Z$.  When $X \to Y$ and $Z \to Y$ are defined relative to some base scheme $Y$, and $f$ is a morphism of $Y$-schemes, then one says $f$ is a \define{projective envelope} if $Z\to Y$ is projective.

\begin{theorem}\label{t.descent}
Let $X \to Y$ be a proper morphism, and let $f\colon Z \to X$ be a projective envelope (relative to $Y$).  Then the sequence
\[
  CK_i(Z \times_X Z) \xrightarrow{{pr_1}_*-{pr_2}_*} CK_i(Z) \xrightarrow{f_*} CK_i(X) \to 0 
\]
is exact, for all $i$.
\end{theorem}

By an application of Chow's lemma, any scheme $X$ admits a projective envelope if it is proper over $Y$.  Further properties of envelopes can be found in \cite{fulton-gillet} or \cite[\S18.3]{fulton}.

Theorem~\ref{t.descent} is proved by following Gillet's argument closely.  In brief, using the terminology of \cite{gillet}, the theorem follows by applying the same arguments to a projective hyperenvelope $Z_\bullet \to X$, using functoriality of the Bousfield-Kan spectral sequence for $K_*(\catM_i(Z_\bullet))$.  We will give an outline here, stressing the novel points for our situtation; see \cite[Appendix~A]{ap} for a similar outline, including most of the terminology and further references.

We will work with an augmented simplicial scheme $Z_\bullet \to X$.  For each non-increasing map of ordinals $\tau\colon \bm \to \bn$ there is structure map $\tau\colon Z_n \to Z_m$.  We require the following condition on sheaves of $\OO_{Z_n}$-modules:
\begin{enumerate}[($*$)]
\item For all $\tau\colon \bm \to \bn$ and all $p>0$, we have $R^p\tau_*\shfF=0$ as sheaves on $Z_m$. \label{g-cond}
\end{enumerate}

Let $\catA_i(Z_n) \subseteq \catM_i(Z_n)$ be the full subcategory of sheaves satisfying \ref{g-cond}.  These are exact categories, and they fit together to form a simplicial category $\catA_i(Z_\bullet)$.  We define
\begin{align*}
  \catE_{q,i}(Z_\bullet) &:= K_q( \catA_i(Z_\bullet) ), 
\end{align*}
where for any category $\catC$, the K-group is defined by Quillen's construction, so $K_q(\catC) = \pi_q( \Omega|NQ\catC| )$.  The natural inclusions of categories $\catA_i(Z_\bullet) \to \catA_{i+1}(Z_\bullet)$ induce homomorphisms $\catE_{q,i}(Z_\bullet) \to \catE_{q,i+1}(Z_\bullet)$.

In our context, the Bousfield-Kan spectral sequence relates $\catE_{i,q}(Z_\bullet)$ with the K-groups $K_q( \catM_i(Z_p) )$.  (Here it is important that each $Z_i$ is quasi-projective over the base $Y$.)  Specifically, it gives a convergent spectral sequence
\[
  E^1_{p,q}(i) = K_q( \catM_i(Z_p) ) \Rightarrow \catE_{p+q,i}(Z_\bullet).
\]
The differential is given by the alternating sum of face homomorphisms.

The key ingredient in proving Theorem~\ref{t.descent} is an analogue of another result of Gillet.

\begin{lemma}\label{l.gillet}
Let $f\colon Z_\bullet \to X$ be a projective hyperenvelope.  Then
\[
  f_*\colon \catE_{q,i} \to K_q(\catM_i(X))
\]
is an isomorphism, natural with respect to the inclusions from $i$ to $i+1$.
\end{lemma}

The proof of the lemma is exactly the same as in \cite{gillet}, everywhere replacing the category of coherent sheaves by the subcategory $\catM_i$.

To deduce Theorem~\ref{t.descent}, one examines the edge homomorphism of the Bousfield-Kan spectral sequence, just as in \cite{gillet}.  Let $Z_1=Z\times_X Z$ and $Z_0 = Z$.  We have a diagram
\[
\begin{tikzcd}
  K_0( \catM_i(Z_1) ) \ar[r] \ar[d] &  K_0( \catM_i(Z_0) ) \ar[r] \ar[d] & E^2_{0,0}(i) \ar[r] \ar[d] & 0 \\
  K_0( \catM_{i+1}(Z_1) ) \ar[r]  &  K_0( \catM_{i+1}(Z_0) ) \ar[r] &  E^2_{0,0}(i) \ar[r]  & 0  
\end{tikzcd}
\]
with exact rows.  From the convergence of the spectral sequence we know $E^2_{0,0}(i)=\catE_{0,i}(Z_\bullet)$, and from Lemma~\ref{l.gillet}, we have $\catE_{0,i}(Z_\bullet) = K_0(\catM_i(X))$.  Putting all this together, the images of the vertical arrows form an exact sequence
\[
  CK_i(Z_1) \to CK_i(Z_0) \to CK_i(X) \to 0,
\]
as claimed.  \qed

\medskip
The descent sequence allows one to extend results from the quasi-projective to the general case.  For instance, suppose $V$ is vector bundle of rank $n$ on a scheme $X$.  The {\it projective bundle formula} asserts that there is an isomorphism of $\ZZ[\beta]$-modules
\[
  CK_*(X) \otimes_{\ZZ[\beta]} CK_*(\PP^{n-1}) \to CK_*(\PP(V)).
\]
In \cite{cai,dl}, this is proved for quasi-projective schemes.  For arbitrary $X$, one can choose a projective envelope $f\colon Z \to X$ (relative to an appropriate base scheme $Y$), and one has a commuting diagram
{\small
\[
\begin{tikzcd}
  CK_*(Z_1) \otimes CK_*(\PP^{n-1}) \ar[r] \ar[d] &  CK_*(Z_0) \otimes CK_*(\PP^{n-1}) \ar[r] \ar[d] &  CK_*(X) \otimes CK_*(\PP^{n-1}) \ar[r] \ar[d] & 0 \\
  CK_*(\PP(V_{1})) \ar[r]  &  CK_*(\PP(V_{0})) \ar[r]  &  CK_*(\PP(V)) \ar[r]  & 0,
\end{tikzcd}
\]
}

\noindent
where $V_0$ and $V_1$ are the pullbacks of $V$ to $Z_0=Z$ and $Z_1=Z\times_X Z$, respectively.  Envelopes are preserved under pullback, so $\PP(V_0) \to \PP(V)$ is an envelope, and therefore the rows are exact, by Theorem~\ref{t.descent}.  The left two vertical arrows are isomorphisms, by the quasi-projective case, so the right vertical arrow is also an isomorphism, by the five lemma.

In particular, this allows one to extend the theory of Chern classes developed in \cite[\S6.7]{cai} to all schemes, as sketched in \cite[Appendix~A]{a}.



\end{document}